\documentclass[11pt]{article}
\usepackage{amsmath, amssymb}
\usepackage{amsthm}
\usepackage{graphicx}
\theoremstyle{definition}
\newtheorem{teo}{Theorem}[section]
\newtheorem{cor}[teo]{Corollary}
\newtheorem{lem}[teo]{Lemma}

\newtheorem{rem}[teo]{Remark}
\newtheorem{prop}[teo]{Proposition}
\newtheorem{defi}[teo]{Definition}

\newcommand{\cl}{\mbox{cl}} 
\newcommand{\conv}{\mbox{conv}} 
\usepackage{enumerate}
\renewcommand{\subset}{\subseteq}
\newcommand{\R}{\mathbb{R}}
\newcommand{\C}{\mathcal{C}}
\newcommand{\D}{\mathcal{D}}
\begin{document}

\title{An incidence Hopf Algebra of Convex Geometries}
\author{Fabi\'an Latorre}
\date{}
\maketitle

A lattice L is ``meet-distributive'' if for each element of L, the meets of the elements directly below it form a Boolean lattice. These objects are in bijection with ``convex geometries'', which are an abstract model of convexity. Do they give rise to an incidence Hopf algebra of convex geometries?

\section{Preliminaries}
We will define some basic concepts in order to understand the question posed.

\begin{defi} A finite lattice $L$ is called \emph{meet distributive} if for every $x\in L-\{\hat 0\}$, the interval $[i(x),x]$ is a boolean lattice, where $i(x)$ is the meet of all elements  $y \in L$ covered by $x$ (i.e. $y<x$ and there is no $z$ such that $y<x<z$). 
\end{defi}

\begin{defi} A closure operator on a set $S$ is a map of subsets of $S$,  $\cl: P(S) \rightarrow P(S)$ such that for any $X,Y\subset S$
\begin{enumerate}[i)]
\item $ X \subset \cl(X)$
\item $ X \subset Y \Rightarrow \cl(X)\subset \cl(Y)$
\item $\cl(\cl(X))=\cl(X)$
\end{enumerate}

\end{defi}

\begin{defi} A closure operator $\cl$ on a set $S$ has the \emph{antiexchange property} if $x \neq y$, $y \in \cl(A \cup x )$, $y \not \in \cl(A) $ then $x \not \in \cl(A \cup y) $.
\end{defi}

\begin{defi}
A convex geometry is a pair $(Z, \cl)$  where $Z$ is a finite set and $\cl$ is a closure operator on $Z$ with the antiexchange property.
\end{defi}

\begin{defi}
On a convex geometry $Z$, $X \subset Z$ is called closed if $\cl(X)=X$. Closed sets of a convex geometry ordered by containment form a lattice $L(Z, \cl)$.
\end{defi}

We have an equivalent definition of convex geometry
\begin{defi} Let $Z$ be a finite nonempty set and $E \subset 2^Z$ a family of subsets of $Z$, then $E$ is called a convex geometry if
\begin{enumerate}[i)]
\item $\emptyset \in E, Z \in E$
\item $X,Y \in E \Rightarrow X \cap Y \in E $
\item $X \in E\backslash \{Z\} \Rightarrow \exists z \in Z \backslash X: X \cup \{z \} \in E$
\end{enumerate}

\end{defi}

The elements of $E$ are called closed subsets of $Z$ and they induce a closure operator on $Z$ where $$\cl(A)= \bigcap_{ \substack{X \in E \\ A \subset X}} X $$
This closure operator has the antiexchange property. Conversely, given $(Z, \cl)$ with the antiexchange property, the family of closed subsets form a family $E$ with the properties defined previously and the closure operator defined as before coincides with $\cl$.\\

The following theorem is due to Edelman \cite{Edelman}. 
\begin{teo}\label{meetconvex}
A finite lattice is meet distributive if and only if it is isomorphic to the lattice of closed sets of a convex geometry.
\end{teo}

In order to construct an incidence Hopf Algebra of meet distributive lattices we must show that the family of such lattices is a hereditary family.

\begin{prop} The family of finite meet distributive lattices is closed under taking subintervals.
\begin{proof} Let $L$ be a meet distributive lattice and $I=[a,b] \subset L$. $I$ is clearly a lattice. Let $x \in I-\{a\} $. There are fewer or an equal number of elements covered by $x$ in $I$ than in $L$ hence their meet $i_I(x) \geq i(x)$, the meet of all elements covered by $x$ in L, so $[i_I(x),x] \subset [i(x),x]$ which is a boolean lattice and subintervals of boolean lattices are boolean lattices. 

\end{proof}

\end{prop}

\begin{prop} The family of meet distributive lattices is closed under direct product.
\begin{proof} Let $L_1,L_2$ be meet distributive lattices. Let $(x_1,x_2) \in L_1 \times L_2 -\{\hat0_1, \hat 0_2\}$. Clearly $i(x_1,x_2)=(i(x_1),i(x_2))$ and $[(i(x_1),i(x_2)), (x_1,x_2)= [i(x_1),x_1]\times [i(x_2),x_2] $ which are boolean lattices and the product of boolean lattices is boolean.

\end{proof}

\end{prop}

\begin{cor} The family of meet distributive lattices is a hereditary family. We denote such family as $\D$. Under isomorphism it forms the incidence Hopf Algebra of meet distributive lattices, which we denote as $M$.

\end{cor}

We call two convex geometries $(Z_1, \cl_1)$, $(Z_2, \cl_2)$ isomorphic if they have isomorphic lattices $L(E_1, \cl_1), L(E_2, \cl_2)$. Then Edelman's theorem gives rise to a bijection between isomorphism classes of meet-distributive lattices and isomorphism classes of convex geometries.
\\

\section{Some examples}
{\bf Example: (Convex Shelling)} Let $P$ be a finite set of points in $\R^n$. Define
$$
E:=\{ X \subset P : \conv(X) \cap P = X \}
$$ 

\begin{figure}[h!]\label{points}
\begin{center}
\includegraphics[scale=0.5]{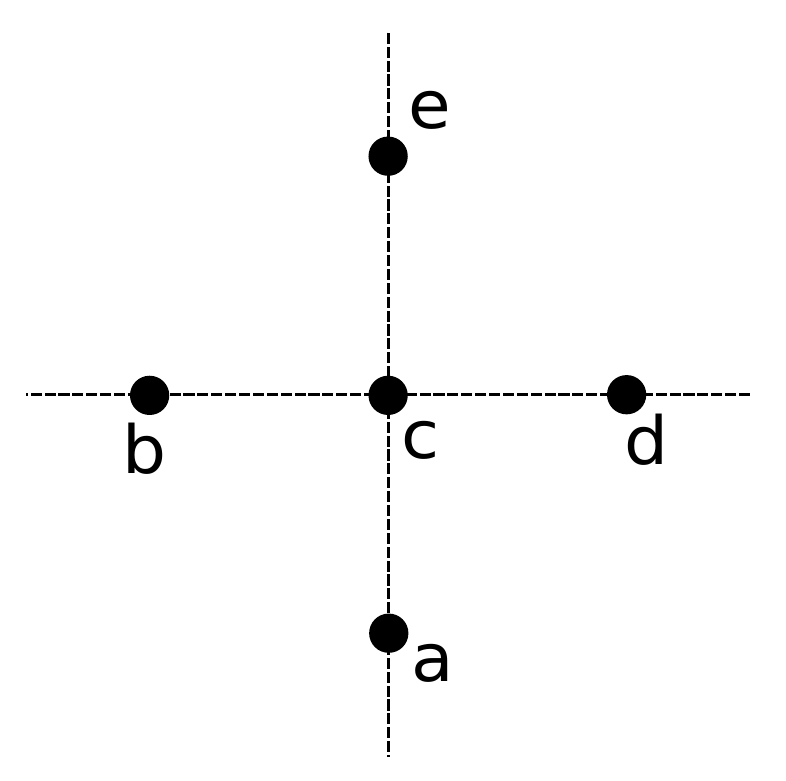} 
\caption{The set of points $P_1$}
\end{center}
\end{figure} 

Where $\conv(X)$ denotes the convex hull of $X$. Then $(P,E)$ is a convex geometry. For example let $P_1=\{a,b,c,d,e\} \subset \R^2$ where $a=(0,-1)$, $b=(-1,0)$, $c=(0,0)$, $d=(1,0)$, $d=(0,1)$ (Figure (\ref{points})). Then the closed sets of $P_1$ in the convex shelling, ordered by size, are $[abcde], [abcd], [abce], [acde], [bcde], [abc]$, $ [acd], [ace]$, $ [bcd], [bce], [cde], [ab], [ac], [ad], [bc], [be], [cd], [ce], [de], [a], [b], [c], [d]$, and the corresponding meet distributive lattice, with the elements in this order, from left to right, and ranked by size is given in figure (2) .\\

\begin{figure}[h!]\label{latticeconv}
\begin{center}
\includegraphics[scale=0.7]{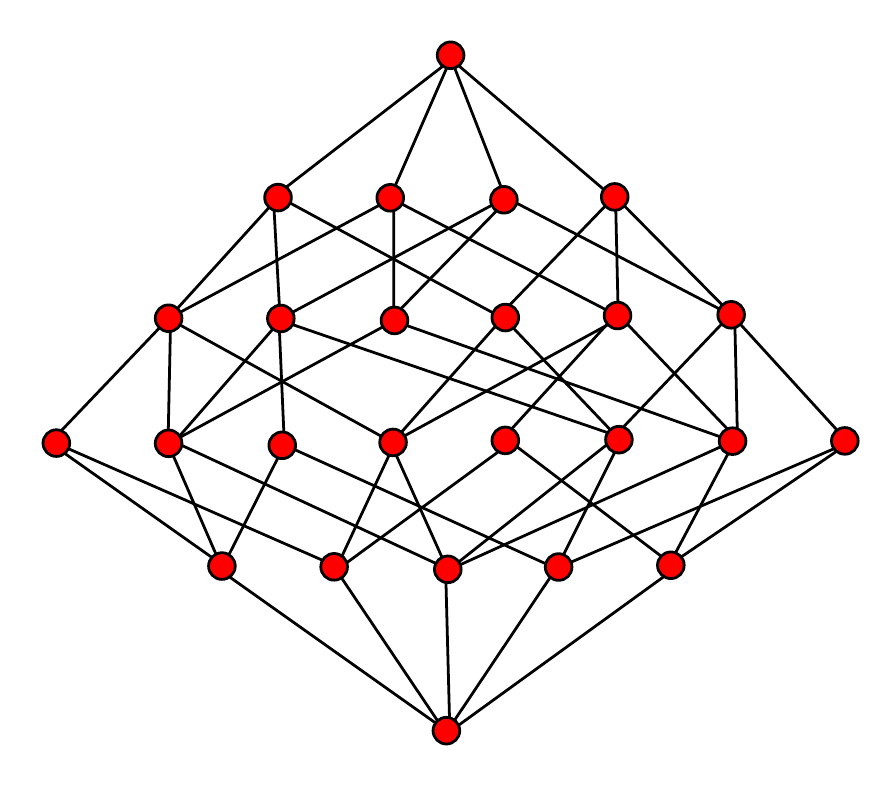} 
\caption{The lattice of the convex shelling $P_1$}
\end{center}
\end{figure} 

{\bf Example (Poset shelling):} Let $P=(E, \leq)$ be a finite poset. Define
$$
E=\{ X \subset E : e \in X, f \leq e \Rightarrow f \in X \}
$$ 
This is the set of downsets of $P$. This is a convex geometry. For example take the power set of $\{1,2\}$ ordered by inclusion. Then the Poset shelling convex geometry is described by the meet distributive lattice in figure (3).

\begin{figure}[h!]\label{poset}
\begin{center}
\includegraphics[scale=0.3]{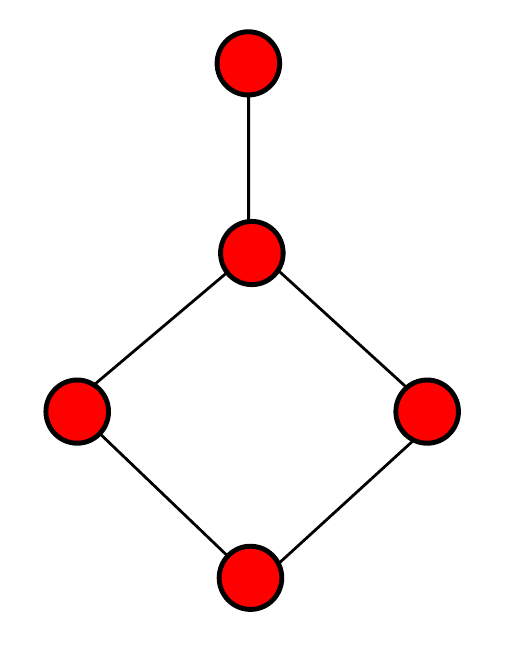} 
\caption{Lattice of down sets of $\{1,2\}$}
\end{center}
\end{figure} 

\section{Convex geometries as a Hopf Algebra}

Now we recall the basic structure of an incidence Hopf algebra. Let $K$ be a field. $K(\D)$ is the free $K$-vector space over $\D$, then the maps
$$
\Delta(L)=\sum_{x\in L} [0_L,x] \otimes [x, 1_L]
$$

$$
\epsilon(L)=\left \{ \begin{array}{cc} 1 & \mbox{if } |L|=1 \\ 0 & \mbox{else} \end{array} \right .
$$

$$
m(L_1,L_2)=L_1 \times L_2
$$

$$
S(L)=\sum_{k \geq 0} \sum_{0=x_0<\ldots <x_k=1} (-1)^k [x_0, x_1]\times \ldots \times [x_{k-1},x_k]
$$

turn $K(\D)$ into a Hopf algebra. We define the following

\begin{defi} Let $(Z,E)$ be a convex geometry over a finite set. and let $A,B$ be closed subsets such that $ A\subset B$. The minor of $(Z,E)$ with respect to $A,B$ is
$$
M[A,B]=\{ X \backslash A : X \in L, A \subset X \subset B \}
$$

\end{defi}

\begin{rem} A minor is always a convex geometry
\end{rem}

\begin{lem}\label{minor} Under Edelman's bijection, any minor $M(A,B)$ corresponds to  the lattice $[A,B] \subset L(Z,E)$, the interval from $A$ to $B$ in the lattice of $(Z,E)$.
\begin{proof}
The map $\phi : M(A,B) \rightarrow [A,B] $, $X \mapsto X \cup A$ is clearly a bijection. If $X,Y \in M(A,B)$, $X \subset Y$, then $X\cup A \subset Y \cup A$ so $\phi(X) \leq \phi(Y)$. \\
\end{proof} 

\end{lem}

\begin{defi} Denote the free $K$-vector space over the set of isomorphism classes of convex geometries as $\C$.

\end{defi}

\begin{cor} The coproduct $\Delta: \C \rightarrow \C \otimes \C$, 
$$(Z,E) \mapsto \sum_{X \in E} M(\emptyset, X)\otimes M(X, Z) $$
 together with the count $\epsilon(Z,E)=1 $ if $Z=\emptyset$, $\epsilon(Z,E)=0$ otherwise, turn $\C$ into a coalgebra.

\end{cor}

\begin{defi} Let $(Z_1, E_1)$ and $(Z_2, E_2)$ be two convex geometries. Suppose $Z_1$ and $Z_2$ are disjoint. Consider the family of subsets of $Z_1 \cup Z_2$, the disjoint union of $Z_1,Z_2$ defined by
$$E_{Z_1Z_2}=\{ X_1 \cup X_2 : X_1 \in E_1, X_2 \in E_2 \}$$. 
\end{defi}

\begin{lem} $(Z_1\cup Z_2, E_{Z_1Z_2}) $ is a convex geometry, called the product convex geometry, for disjoint $Z_1$, $Z_2$.
\begin{proof}

\begin{enumerate}[i)]
\item $\emptyset \in E_1, \emptyset \in E_2 \Rightarrow \emptyset=\emptyset \cup \emptyset \in E_{Z_1Z_2}$. $Z_1 \in E_1, Z_2 \in E_2 \Rightarrow Z_1\cup Z_2 \in E_{Z_1Z_2}$
\item Suppose $X_1\cup X_2, Y_1\cup Y_2 \in E_{Z_1,Z_2} $ then $X_1 , Y_1 \in E_1$ and $X_1 \cap Y_1 \in E_1$ as $(Z_1,E_1)$ is a convex geometry. Analogously $X_2 \cap Y_2 \in E_2$ then $(X_1\cup X_2) \cap (Y_1\cup Y_2)= (X_1 \cap Y_1) \cup (X_2 \cap Y_2) \in E_{Z_1Z_2}$. 
\item Suppose $X_1\cup X_2 \in E_{Z_1Z_2}\backslash \{Z_1\cup Z_2 \} $ then either $X_1 \in E_1\backslash \{Z_1\} $ or  $X_2 \in E_1\backslash \{Z_2\} $.  Suppose $X_1 \in E_1\backslash \{Z_1\} $, then there exists some $z \in Z_1$ such that $X_1\cup \{z\} \in E_1$  so that $X_1 \cup X_2 \cup \{z\}=(X_1\cup\{z\})\cup X_2 \in E_{Z_1Z_2}$. The other case follows by symmetry. \\
\end{enumerate}

\end{proof}

\end{lem} 

\begin{lem}\label{product} For disjoint convex geometries $(Z_1, E_1), (Z_2, E_2)$ the lattice $L(Z_1 \cup Z_2)$ is isomorphic to $L(Z_1)\times L(Z_2)$.
\begin{proof}
Define $\phi : L(Z_1) \times L(Z_2) \rightarrow L(Z_1 \cup Z_2)$ as
$$
\phi(X,Y)=X\cup Y
$$
This map is clearly a bijection. Now if $ (X_1,Y_1) \leq (X_2, Y_2)$ in $ L(Z_1)\times L(Z_2)$, then $X_1\subset X_2$ and $Y_1 \subset Y_2$. It follows that
$$
\phi(X_1,Y_1)=X_1\cup Y_1 \subset X_2 \cup Y_2 = \phi(X_2,Y_2)
$$

\end{proof}

\end{lem}

\begin{cor} The multiplication map $m: \C \otimes \C \rightarrow \C $ given by
$$
m(Z_1 \otimes Z_2) = Z_1 \cup Z_2
$$
Turns $\C$ into a bialgebra.
\end{cor}

\begin{rem} Although dealing with isomorphism classes, this multiplication map is well defined by the previous lemma and the fact that one can choose representatives with disjoint $Z_1$ and $Z_2$.

\end{rem}

\section{Some sub-Hopf Algebras and further considerations}

It is a matter of interest if certain subclasses of convex geometries form a sub-Hopf algebra of the general incidence Hopf algebra of convex geometries. We cite the following theorems

\begin{teo}\label{teo1} {\bf Nakamura-Okamoto (2003)} The convex geometry class of convex shellings of finite point sets is not closed under taking minors.

\end{teo} 

\begin{teo} {\bf Nakamura-Okamoto (2003)}\label{teo2} The class of digraph point search convex geometries is closed under taking minors

\end{teo}

\begin{teo} {\bf Nakamura (2003)}\label{teo3}The class of poset shellings convex geometries is closed under taking minors

\end{teo}

\begin{teo} {\bf Nakamura (2003) }\label{teo4} $(Z,E)$ is a poset shelling if and only if it contains no minor isomorphic to the following:
\begin{figure}[h!]\label{poset}
\begin{center}
\includegraphics[scale=0.3]{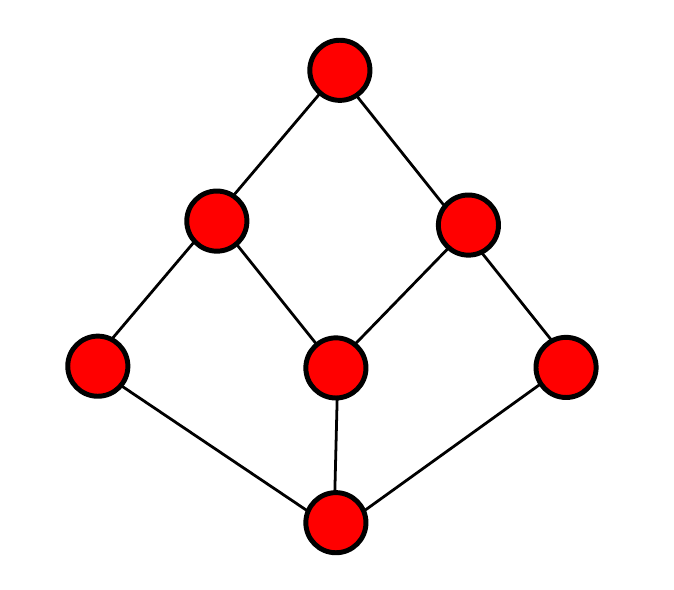} 
\caption{Forbidden minor for poset shelling convex geometries $F_1$}
\end{center}
\end{figure} 

\end{teo}

\begin{teo}\label{teo5} {\bf Nakamura-Okamoto} $(Z,E)$ is a digraph point search convex geometry if and only if it contains no minor isomorphic to the following:

\begin{figure}[h!]\label{poset}
\begin{center}
\includegraphics[scale=0.3]{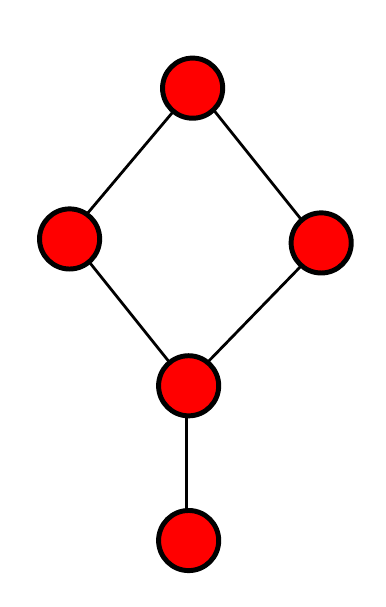} 
\caption{Forbidden minor for digraph point search convex geometries $F_2$}
\end{center}
\end{figure} 

\end{teo}

\begin{cor} By theorem (\ref{teo1}), the class of convex shellings does not yield a sub-Hopf algebra of $\C$.

\end{cor}

Now we prove that both the poset shellings and the digraph point search classes are closed under products.

\begin{teo} The Convex geometry class of poset shellings and the class of digraph point search are closed under products.

\begin{proof} Let $(Z_1, E_1), (Z_2,E_2)$ be two convex geometries of the same class. Let $L_1$ and $L_2$ be their corresponding closed subsets lattice. By lemma (\label{minor}), the minors of such convex geometries correspond to subintervals in the corresponding lattice. By theorems (\ref{teo4}),(\ref{teo5}) we then know that both $L_1$ and $L_2$ contain no subinterval isomorphic to either $F_1$, in the case of poset shellings, or $F_2$, in the case of digraph point search. Now by lemma (\ref{product}) we know the lattice of closed sets of $Z_1\cup Z_2$ corresponds to the product of the lattices $L_1\times L_2$. \\

Any subinterval $I$ of $L_1 \times L_2$ is isomorphic to a product $I_1 \times I_2$ where $I_1$ is a subinterval of $L_1$ and $I_2$ is a subinterval of $L_2$, Also we have the equality $|I|=|I_1||I_2|$. Now Suppose $L_1 \times L_2$ contains an interval $I$ isomorphic to either $F_1$ or $F_2$ then we have a decomposition $|F_i|=|I_1||I_2|$, but the size of $F_i$ is prime in both cases so either $|I_1|=1$ or $|I_2|=1$. Without loss of generality $|I_2|=1$. As $I_2$ is a one element lattice we have that $F_i\simeq I_1 \times I_2 \simeq I_1 $ which is a contradiction as $L_1$ contains no subinterval isomorphic to $F_i$. \\

\end{proof} 

\end{teo}

\begin{cor} The class of poset shellings and the class of digraph point search convex geometries form two sub-Hopf Algebras of the general incidence Hopf Algebra of convex geometries $\C$. 

\end{cor}

Now consider the set $Z_n:=[n]$ and the family of subsets $E=\{ [k] : 1 \leq k \leq n \} \cup \{ \emptyset \} $. It is straightforward to check that $(Z,E)$ is a convex geometry, moreover, it forms a finite linear order and any finite linear order is isomorphic to some $Z_n$. Therefore we also have

\begin{lem} The Hopf algebra of symmetric functions is a sub Hopf algebra of $\C$.
\begin{proof}
Consider the map
$$
\phi( \mbox{finite linear order of size $n$} ) = Z_{n-1}
$$
This is extends to an injective Hopf algebra map from $Sym$ to $\C$.
\end{proof}

\end{lem}

The set of isomorphism classes of convex geometries is really large. We don't expect a simplified formula for the antipode, however we may find simplified formulas for the antipode on certain families of convex geometries.

{\bf Example } Consider the family of convex shellings of regular $n$-gons $\{ P_n \} $ in two dimensional space. We call the convex hull of two different points the $2$-gon and a single point the $1$-gon. For the regular $n$-gon, we see that any subset of the vertices belongs to the convex shelling, so we can distinguish an individual vertex $v$ and separate the elements of the convex shelling as those containing $v$ and those not containing $v$. The elements of the convex shelling that do not contain $v$ is the convex shelling of a $(n-1)$-gon. \\

We see at once that the convex shelling of the $n$-gon is the same as the product of the convex shelling of the $(n-1)$-gon and the convex shelling of a vertex, $P_1$. Hence by induction we have that $P_n= (P_1)^n$. Then $S(P_n)=S(P_1)^n$.

Now $P_1$ corresponds to the lattice of a linear order of size $2$, $L_1$. Using the formula of the antipode we have $S(P_1)=S(L_1)=-L_1=-P_1$ then 

$$
S(P_n)=(-1)^n P_n
$$ 

An analogous argument with some other considerations would yield a formula for the antipode of the convex shellings of the vertices of simplicial polytopes (every facet is a simplex) on higher dimension.


\newpage
\begin{thebibliography}{100}

\bibitem{Edelman} Edelman, Paul H, (1980). Meet-distributive lattices and the anti-exchange closure, Algebra Universalis 10 (1): 290Ð299.


\bibitem{okamoto} Yoshio Okamoto, Masataka Nakamura, (2003). The forbidden minor characterization of line search antimatroids of rooted digraphs. Discrete Applied Mathematics 131(2): 523-533

\bibitem{nakamura} Masataka Nakamura, (2003). Excluded-minor characterizations of antimatroids arisen from posets and graph searches. Discrete Applied Mathematics 129 (2Ð3): 487Ð498 

\end{thebibliography}
\end{document}